\magnification=1200

\def\C{{\bf C}}

\def\N{{\bf N}}
\def\R{{\bf R}}

\def\hal{{\vrule height 10pt width 4pt depth 0pt}}
\def\la{{\langle}}
\def\ra{{\rangle}}

\centerline{\bf Hilbert bimodules with involution}
\medskip

\centerline{Nik Weaver}
\bigskip
\bigskip

{\narrower{
\noindent {\it We examine Hilbert bimodules which possess a (generally
unbounded) involution. Topics considered include a linking algebra
representation, duality, locality, and the role of these bimodules in
noncommutative differential geometry.}
\bigskip}}
\bigskip

Many Hilbert modules which arise in practice carry natural involutions,
typically deriving from the involutions of C*-algebras involved in their
construction. Usually these Hilbert module involutions are not only
non-isometric, they are unbounded --- possibly even if the module is
finitely generated.
\medskip

The same examples generally also have a bimodule structure which interacts
with, and may be recovered from, the involution via the equation $(ax)^* =
x^*a^*$, where $a$ is an element of the C*-algebra and $x$ is an element of
the Hilbert module. Thus, involutions are closely related to bimodule
structure. (But these are very different from the sort of bimodules that
arise in the context of Morita equivalence [15].)
\medskip

Philosophically, if one regards Hilbert modules as ``noncommutative complex
Hilbert bundles'' ([16], [20], [21]) then Hilbert bimodules with involution
may be seen as ``noncommutative real Hilbert bundles.'' Indeed, assuming a
locality condition, in the commutative case the extra structure provided by
the involution corresponds precisely to a real structure on the corresponding
bundle (Theorem 11).
\medskip

The motivating commutative example is the tangent bundle $TX$ of a Riemannian
manifold $X$. Since the tangent space at each point carries an inner product,
the space $S_0(TX)$ of continuous sections of $TX$ which vanish at infinity
has a $C_0(X)$-valued inner product given by $\langle \phi, \psi\rangle(x) =
\langle \phi(x), \psi(x)\rangle$. If we are using complex scalars, we
must complexify $S_0(TX)$ to make it a module over $C_0(X)$, and we must
also extend the inner product to the complexification. But in addition the
complexified module now has an involution given by $\phi_0 + i\phi_1 \mapsto
\phi_0 - i\phi_1$. Thus, in the noncommutative setting ([4], [17], [23])
we expect a description of ``noncommutative Riemannian structure'' to
involve involutive Hilbert bimodules.
\medskip

This description should also include an ``exterior derivative'' realized as
an unbounded self-adjoint derivation from the C*-algebra into the bimodule.
Therefore both the $*$-operation and the bimodule structure are important.
\medskip

It was argued in [13] that operator modules are the correct noncommutative
version of complex Banach bundles. Thus self-adjoint operator bimodules
may be seen as a noncommutative version of real Banach bundles. But we do
not pursue this issue here.
\medskip

The proof of Theorem 8 was supplied by Charles Akemann, using the powerful
excision technique of [2].
\bigskip
\bigskip

\noindent {\bf 1. Preliminaries.}
\bigskip

We will use basic facts about Hilbert modules without comment; see [10], [11],
or [14] for background.
\medskip

In the following definition we require the involution to be defined
everywhere. Thus, we get around the unboundedness problem mentioned in
the introduction by eliminating that portion of the module on which the
involution is not defined. This means that what we call the left and right
seminorms, $\|\cdot\|_l$ and $\|\cdot\|_r$, in general cannot
be complete. However, there is no obstruction to completeness of the max
norm $\|\cdot\|_m$ and in general this seems to be the appropriate requirement.
\bigskip

\noindent {\bf Definition 1.} Let $A$ be a pre-C*-algebra. A {\it pre-Hilbert
$*$-bimodule over $A$} is an $A$-$A$-bimodule $E$ together with an $A$-bilinear
map $\la \cdot, \cdot \ra: E \times E \to A$ and an antilinear involutive map
${}^*: E \to E$ such that
\medskip

(a) $\la x, y \ra^* = \la y^*, x^* \ra$,
\medskip

(b) $(ax)^* = x^* a^*$, and
\medskip

(c) $\la x, x^* \ra \geq 0$
\medskip

\noindent hold for all $a \in A$ and $x,y \in E$. We define sesquilinear
$A$-valued inner
products $\la x, y \ra_l = \la x, y^* \ra$ and $\la x, y \ra_r =
\la x^*, y \ra$ and seminorms $\|x\|_l^2 = \| \la x, x \ra_l \|$,
$\|x\|_r^2 = \| \la x, x \ra_r \|$, and $\|x\|_m = \max(\|x\|_l, \|x\|_r)$.
If $A$ is a C*-algebra, we say $E$ is a {\it Hilbert $*$-bimodule over $A$}
provided $\|\cdot\|_m$ is a complete norm.\hfill\hal
\bigskip

By $A$-bilinearity we mean that $\la ax, y \ra = a \la x, y \ra$,
$\la xa, y \ra = \la x, ay \ra$, and $\la x, ya \ra = \la x, y \ra a$ for
$a \in A$ and $x,y \in E$.
\medskip

In most cases $A$ will be complete from the start, but in one or two places we
will want to allow it to be incomplete. (We never need to allow nonzero
elements in $A$ to have zero norm, however.) To make sense of axiom (c), we
take positivity in an incomplete algebra to mean positivity in its completion.
\medskip

If $A$ is complete then $E$ is a left pre-Hilbert $A$-module with
respect to the inner product $\la \cdot, \cdot \ra_l$ and a right
pre-Hilbert $A$-module with repect to $\la \cdot, \cdot \ra_r$. In any case
$\|\cdot\|_l$ and $\|\cdot\|_r$ are seminorms by the same argument which
shows this for ordinary Hilbert modules.
\medskip

Note also that $\|x\|_l = \|x^*\|_r$ for any $x \in E$; in particular, if $x$
is self-adjoint then $\|x\|_l = \|x\|_r = \|x\|_m$.
\medskip

Our first order of business is to show that pre-Hilbert $*$-bimodules can
always be completed to Hilbert $*$-bimodules.
\bigskip

\noindent {\bf Lemma 2.} {\it Let $A$ be a pre-C*-algebra and $E$ a
pre-Hilbert $*$-bimodule over $A$. Then for any $a \in A$ and $x \in E$ we have
$$ \|ax\|_l, \|xa\|_l \leq \|a\| \|x\|_l\qquad{\rm and}\qquad
\|ax\|_r, \|xa\|_r \leq \|a\| \|x\|_r.$$}

\noindent {\it Proof.} First, we have
$$\|ax\|_l^2 = \| \la ax, (ax)^* \ra \| = \|a \la x, x^* \ra a^*\|
\leq \|a\|^2 \|x\|_l^2,$$
so $\|ax\|_l \leq \|a\| \|x\|_l$.
To see that $\|xa\|_l \leq \|a\| \|x\|_l$, observe that $b \leq c$ implies
$\la xb, x \ra_l \leq \la xc, x \ra_l$ since
$$\la x(c-b), x \ra_l = \la x(c - b)^{1/2}, x(c - b)^{1/2} \ra_l \geq 0.$$
Without loss of generality suppose $\|a\| \leq 1$; then letting $b = aa^*$
we have $b^2 \leq b$ and
$$\eqalign{0 &\leq \la x - xb, x - xb \ra_l\cr
&= \la x, x \ra_l - 2\la xb, x \ra_l + \la xb^2, x \ra_l\cr
&\leq \la x, x \ra_l - \la xb, x \ra_l\cr}$$
so that $\|xa\|_l^2 = \| \la xb, x \ra_l \| \leq \|x\|_l^2$ as desired.
Taking adjoints yields the same inequalities for the norm
$\|\cdot\|_r$.\hfill\hal
\bigskip

\noindent {\bf Proposition 3.} {\it Let $E$ be a pre-Hilbert $*$-bimodule
over a pre-C*-algebra $A$ and let $N = \{x \in E: \|x\|_m = 0\}$.
Then $N$ is a sub-bimodule of $E$ and the inner product and involution
on $E$ descend to $E/N$ and extend to the completion of $E/N$. The
completion of $E/N$ is a Hilbert $*$-bimodule over the completion of $A$.}
\medskip

\noindent {\it Proof.} By Lemma 2, $N$ is a sub-bimodule of $E$
and the left and right actions of $A$ on $E/N$ extend continuously to its
completion for $\|\cdot\|_m$; applying Lemma 2 again allows us to extend the
module actions to the completion of $A$.
The inner product descends to $E/N$ and then extends to its
completion by the Cauchy-Schwarz inequality for ordinary Hilbert modules,
and the corresponding assertions for the involution are trivial.
Axioms (a) to (c) of Definition 1 all hold in the completion by
continuity.\hfill\hal
\bigskip

Next, we briefly consider bounded module maps on Hilbert $*$-bimodules.
\bigskip

\noindent {\bf Definition 4.} Let $E$ be a Hilbert $*$-bimodule over
a C*-algebra $A$. We let $E_l$ and $E_r$ denote the set $E$ considered
respectively as a left or right pre-Hilbert $A$ module with the left or right
inner product described in Definition 1.
If $F$ is another Hilbert $*$-bimodule over $A$ then we define
$B_l(E, F)$ to be the set of all left $A$-linear adjointable maps from
$E_l$ into $F_l$ such that
$$\|T\| = \sup_{x \in E} \{ \|Tx\|_l/\|x\|_l, \|Tx\|_r/\|x\|_r \} < \infty,$$
assuming the convention $0/0 = 0$.
\medskip

We define $B_r(E, F)$ similarly as the right $A$-linear adjointable maps from
$E_r$ to $F_r$ and write $B_l(E) = B_l(E, E)$ and $B_r(E) =
B_r(E, E)$.\hfill\hal
\bigskip

For $a \in A$ define $L_a, R_a: E \to E$ by $L_a x = ax$ and $R_a x = xa$.
\bigskip

\noindent {\bf Proposition 5.} {\it Let $A$ be a C*-algebra and let $E$ and
$F$ be Hilbert $*$-bimodules over $A$. Then $B_l(E, F)$ is naturally
isometrically anti-isomorphic to $B_r(E, F)$. Furthermore, $B_l(E)$ and
$B_r(E)$ are C*-algebras, and
the left and right representations $a \mapsto L_a$, $a \mapsto R_a$ are
$*$-homomorphisms from $A$ into $B_r(E)$ and $B_l(E)$, respectively.}
\medskip

\noindent {\it Proof.} Given $T: E \to F$ define $\widehat{T}: E \to F$ by
$\widehat{T}(x) = (Tx^*)^*$. If $T \in B_l(E, F)$ then $\widehat{T}(xa) =
(\widehat{T} x)a$, so $\widehat{T}$ is right $A$-linear, and
$$\la \widehat{T} x, y\ra_r = \la T(x^*), y \ra = \la T(x^*), y^* \ra_l
= \la x^*, T^*y^* \ra_l = \la x, \widehat{T^*} y \ra_r,$$
so $\widehat{T}$ is adjointable. Thus the map $T \mapsto \widehat{T}$ is
a linear anti-isomorphism from
$B_l(E, F)$ into $B_r(E, F)$, and since norms are computed the same way in
each case the correspondence is isometric.
\medskip

It is trivial to check that $B_l(E)$ is complete.
Let $B_l^0(E)$ be the pre-C*-algebra of left $A$-linear adjointable maps
bounded only for the seminorm $\|\cdot\|_l$, and define $B_r^0(E)$ similarly
using only the seminorm $\|\cdot\|_r$. (These spaces may not be complete
because $E$ need not be complete for the left and right seminorms separately.)
The map $T \mapsto T \oplus \widehat{T}^*$ isometrically embeds $B_l(E)$ as a
C*-subalgebra of $B_l^0(E) \oplus B_r^0(E)$. Since $B_r(E)$ is the opposite
algebra of $B_l(E)$ it is a C*-algebra too.
\medskip

Lemma 2 shows that $L_a \in B_r(E)$ and $R_a \in B_l(E)$ for all
$a \in A$. Linearity and multiplicativity of the representations are
clear. Preservation of adjoints follows from the calculation
$$\la ax, y \ra_r = \la x^* a^*, y \ra = \la x^*, a^* y \ra
= \la x, a^* y \ra_r.$$
(Hence $L_{a^*} = L_a^*$, and $R_{a^*} = R_a^*$ is proven
analogously.)\hfill\hal
\bigskip
\bigskip

\noindent {\bf 2. Simple examples.}
\bigskip

We now list several simple examples of Hilbert $*$-bimodules.
\bigskip

\noindent {\bf Example 6.} {\it The case $A = \C$.} For $X$ any measure
space, $L^2(X)$ has a canonical bilinear $\C$-valued map $(f,g) \mapsto
\int fg$ and involution $f \mapsto \bar{f}$.
\medskip

If $H$ is any real Hilbert space then the complex Hilbert space $H + iH$
has the bilinear form
$$(v_1 + iw_1, v_2 + iw_2) \mapsto \la v_1, v_2\ra + i\la w_1, v_2\ra
+ i\la v_1, w_2\ra - \la w_1, w_2\ra$$
and carries the natural involution $v + iw \mapsto v - iw$. 
\medskip

Let $H$ be a complex Hilbert space and let $H^*$ be its dual space; for
$v \in H$ write $\delta_v$ for the linear functional $w \mapsto \la w, v\ra$.
Then $H \oplus H^*$ with bilinear form
$$(v_1 \oplus \delta_1, v_2 \oplus \delta_2) = \delta_2(v_1)$$
and involution $v \oplus \delta_w \mapsto w \oplus \delta_v$ provides a
simple example
where the seminorms $\|\cdot\|_l$ and $\|\cdot\|_r$ do not agree. Indeed
$\|v \oplus 0\|_l = \|v\|$ and $\|v \oplus 0\|_r = 0$ for any $v \in H$.
\medskip

If $A$ is any C*-algebra and $\phi$ is a state on $A$, then the bilinear
map $(x,y) \mapsto \phi(xy)$ together with the original involution on $A$
makes $A$ a pre-Hilbert $*$-bimodule over $\C$. This is just the GNS
construction and here the left and right
seminorms coincide only if $\phi$ is a trace. Similarly, if $M$ is a
von Neumann algebra and $\phi$ is a weight on $M$ then the same prescription
makes a pre-Hilbert $*$-bimodule over $\C$ of the set $\{x \in M: \phi(xx^*),
\phi(x^* x) < \infty\}$.\hfill\hal
\bigskip

\noindent {\bf Example 7.} {\it Examples derived from tensor products.}
Let $A$ be a C*-algebra, $I$ a closed ideal of $A$, and $H$ a real Hilbert
space. Then the algebraic tensor product $I \otimes_\R H$ is a pre-Hilbert
$*$-bimodule over $A$. The bilinear map into $A$ is defined by
$(a \otimes v, b \otimes w) \mapsto \la v, w\ra ab$ and the involution by
$(a \otimes v)^* = a^* \otimes v$. Taking $H = l^2(\N)$ and $I = A$ in this
construction yields the uncompleted ``standard'' Hilbert module
$l^2(A)$. Completing recovers precisely the sequences $(a_n) \subset A$
such that $\sum a_n^*a_n$ and $\sum a_na_n^*$ both converge in norm.\hfill\hal
\bigskip

In Example 7, if $I = A$ and $H = \R$ then the left and right seminorms
coincide since $\|xx^*\| = \|x^*x\|$. However, in general they disagree, even
for $H = \R^2$. In fact this happens whenever $A$ is noncommutative, by the
following theorem of Akemann [1].
\bigskip

\noindent {\bf Theorem 8 (Akemann).} {\it Let $A$ be a C*-algebra that is not
commutative. Then there exist $x,y \in A$ such that $\|xx^* + yy^*\|
\neq \|x^*x + y^*y\|$.}
\medskip

\noindent {\it Proof.} Since $A$ is not commutative there exists a pure
state $f$ which gives rise to an irreducible representation on a Hilbert
space of dimension at least two. Let $(a_\alpha)$ be a decreasing net
of positive norm one elements which satisfies $f(a_\alpha) = 1$ for all
$\alpha$ and excises $f$ ([2], Proposition 2.2). This means that
$$\| a_\alpha b a_\alpha - f(b) a_\alpha^2 \| \to 0$$
for every $b \in A$.
\medskip

By Kadison's transitivity theorem ([12], Theorem 2.7.5) we can find a unitary
$u \in A$ such that $f(u^* a_{\alpha_0} u) < 1$ for some fixed $\alpha_0$.
Define $b_\alpha = u^* a_\alpha u$. Then $\| f(b_{\alpha_0})a_\alpha^2 \|
< 1$ and so the excision condition implies that
$\|a_\beta b_{\alpha_0} a_\beta\| < 1$ for some $\beta > \alpha_0$. Since
the net is decreasing, it follows that $\|a_\beta b_\beta a_\beta\| < 1$.
Thus $\|a_\beta b_\beta^{1/2}\|, \|b_\beta^{1/2} a_\beta\| < 1$, hence
$\|a_\beta b_\beta\|, \|b_\beta a_\beta\| < 1$.
\medskip

Set $x = u^* a_\beta$ and $y = a_\beta$. Then
$$\|x^*x + y^*y\| = \|a_\beta^2 + a_\beta^2\| = 2$$
but
$$\eqalign{\|xx^* + yy^*\|
&= \|b_\beta^2 + a_\beta^2\| \leq \|b_\beta + a_\beta\|\cr
&= \|b_\beta^2 + b_\beta a_\beta + a_\beta b_\beta + a_\beta^2\|^{1/2}\cr
&< (1 + 1 + 1 + 1)^{1/2} = 2.\cr}$$
Thus $\|x^*x + y^*y\| \neq \|xx^* + yy^*\|$.\hfill\hal
\bigskip

\noindent {\bf Example 9.} {\it Examples arising from maps into subalgebras.}
If $B \subset A$ are C*-algebras and $\phi: A \to B$ is a conditional
expectation, then $A$ is a $B$-$B$-bimodule and the maps $(x,y) \mapsto
\phi(xy)$ and $x \mapsto x^*$ make it a pre-Hilbert $*$-bimodule.
\medskip

Similarly, if $N \subset M$ are von Neumann algebras and $\phi: M_+ \to
\widehat{N_+}$ is an operator-valued weight [6], then $E = \{x \in M:
\|\phi(x^*x)\|, \|\phi(xx^*)\| < \infty\}$ is a pre-Hilbert $*$-bimodule over
$N$ with bilinear form $(x,y) \mapsto \hat{\phi}(xy)$, where $\hat{\phi}$ is
the unique linear extension of $\phi$.\hfill\hal
\bigskip

\noindent {\bf Example 10.} {\it The local commutative case.} Let $X$ be a
locally compact space and let $B$ be a Fell bundle of real Hilbert spaces
over $X$ [7]. Let $S_0(B)$ be the space of continuous sections of $B$ which
vanish at infinity, and let $S_0^\C(B) = S_0(B) + iS_0(B)$ be its
complexification. Then  $S_0^\C(B)$ is a bimodule over $C_0(X)$ with
coincident left and right actions given by fiberwise multiplication.
The involution on $S_0^\C(B)$ is defined by $(f + ig)^* = f - ig$
and the bilinear form by bilinear extension of the fiberwise inner product
on $S_0(B)$.\hfill\hal
\bigskip

An important special case of Example 10 arises when $X$ is a Riemannian
manifold and $B$ is its tangent bundle, as mentioned in the introduction.
This has a simple noncommutative generalization when we have an action of a
real Lie group $G$ on a C*-algebra $A$ [4]. Fixing an inner product on the
Lie algebra $g$ of $G$, the Hilbert $*$-bimodule $E = A \otimes_\R g$ formed
as in Example 7 plays the role of the tangent bimodule. Sauvageot's
construction generalizes this class of examples (see Section 4).
\medskip

Example 10 has the following converse.
\bigskip

\noindent {\bf Theorem 11.} {\it Let $A = C_0(X)$ be a commutative C*-algebra
and let $E$ be a Hilbert $*$-bimodule over $A$. Suppose that any inner product
of self-adjoint elements of $E$ is a self-adjoint element of $A$. Then there
is a Fell bundle $B$ of real Hilbert spaces over $X$ such that $E \cong
S_0^\C(B)$.}
\medskip

\noindent {\it Proof.} Let $Y$ be the one-point compactification of $X$
and let $A^\sim = C(Y)$. Then $E$ is also a Hilbert $*$-bimodule over $A^\sim$.
We want to show that the left and right actions of
$A^\sim$ on $E$ coincide. Fix $p \in Y$. For any $a \in A^\sim$ such that $a \geq 0$
and $a(p) = 0$, and any $x \in E$ such that $x = x^*$, set $b = a^{1/2}$; then
$$\eqalign{\la bx + xb, bx + xb\ra(p)
&= b(p)\la x, bx\ra(p) + b(p)\la x, xb\ra +
\la xb, bx\ra(p) + \la xb, x\ra(p)b(p)\cr
&= \la xb, bx\ra(p),\cr}$$
so that $\la x, ax\ra(p) \geq 0$. But also
$$\eqalign{\la bx + xb, ibx - ixb\ra(p)
&= ib(p)\la x, bx\ra(p)  - ib(p)\la x, xb\ra(p)
+ i\la xb, bx\ra(p) - i\la xb, x\ra(p)b(p)\cr
&= i\la xb, bx\ra(p).\cr}$$
Since both $bx + xb$ and $ibx - ixb$ are self-adjoint, so is their inner
product, so this computation shows that $\la x, ax\ra(p)$ must be
purely imaginary. But we already showed that it is real, so it follows that
$\la x, ax\ra(p) = 0$. By linearity, we have $\la x, ax\ra(p) = 0$ for any
$a \in A^\sim$ such that $a(p) = 0$.
\medskip

Now let $a \in A^\sim$ and $x \in E$ and suppose $a$ is real and $x = x^*$.
Since the quantity
$$\la ax + xa, iax - ixa\ra = i\la x, a^2x\ra - ia^2\la x,x\ra$$
is real and $a^2\la x,x\ra$ is also real, it follows that
${\rm Re} \la x, a^2x\ra = a^2\la x,x\ra$. Therefore
$$\eqalign{\la ax - xa, ax - xa\ra &= {\rm Re}\la ax - xa, ax - xa\ra\cr
&= a^2\la x,x\ra - 2{\rm Re}\, a\la x, ax\ra  + {\rm Re}\la x, a^2x\ra\cr
&= 2{\rm Re}(a^2\la x,x\ra - a\la x, ax\ra),\cr}$$
and evaluating this expression at $p \in Y$ yields
$$2{\rm Re}(a(p)\la x, (a(p) - a)x\ra(p)),$$
which is zero by the last paragraph. Since this is true for all $p \in Y$,
we have $ax - xa = 0$ as desired. Taking linear combinations, we conclude
that this is true for any $a \in A^\sim$ and $x \in E$.
\medskip

Define $E_{sa} = \{x \in E: x = x^*\}$, so that $E = E_{sa} + iE_{sa}$.
Then $E_{sa}$ is a real Hilbert module over $C(Y; \R)$, hence
$E_{sa} \cong S(B)$ for some real Hilbert Fell bundle over $Y$ [21]. Since the
inner product on $E$ takes values in $C_0(X)$, this bundle must have zero
fiber over the point at infinity, so actually $E_{sa} \cong S_0(B)$ for
some Hilbert Fell bundle $B$ over $X$. Thus $E \cong S_0^\C(B)$.\hfill\hal
\bigskip

The left and right actions do not coincide in general in the commutative
case; for instance consider a module constructed as in Example 9 using a
conditional expectation from a C*-algebra onto a commutative but not
central subalgebra. Also, the case of $H\oplus H^*$ in Example 6 shows that
the reality condition on inner products can fail even if the left and
right actions agree. However, the construction in Example 10 always
verifies the reality condition, so Theorem 11 is a true converse of it.
\bigskip
\bigskip

\noindent {\bf 3. Operator representation.}
\bigskip

Let $A$ be a C*-algebra and $E$ a Hilbert $*$-bimodule over $A$. The following
is the appropriate version of the linking algebra construction in this setting.
Let $N = \{x \in E: \|x\|_r = 0\}$ and let $E_r$ be the completion of $E/N$
for $\|\cdot\|_r$, so that $E_r$ is a right
Hilbert module over $A$. The $A$-valued inner product on $E_r$ extends
$\langle \cdot, \cdot\rangle_r$, and we use the same notation for the
extension. Then define $F = A \oplus E$ to be the direct sum
of right Hilbert $A$-modules.
\medskip

Let $B(F)$ be the space of bounded adjointable right $A$-linear maps from
$F$ to itself. This is a C*-algebra. Then define $\phi: A \to B(F)$ by
$$\phi(a)(b \oplus y) = ab \oplus ay$$
and $\psi: E \to B(F)$ by
$$\psi(x)(b \oplus y) = \langle x^*,y\rangle_r \oplus xb.$$
A number of simple facts need to be verified. First,
for any $a \in A$ the map $\phi(a)$ is bounded by Lemma 2, and a short
computation shows that $\phi(a)^* = \phi(a^*)$, so $\phi$
is an injective $*$-homomorphism from $A$ into $B(F)$. For any $x \in E$ the
map $\psi(x)$ is clearly linear, and it is bounded because
$$\eqalign{\langle \psi(x)(b\oplus y), \psi(x)(b\oplus y)\rangle_r
&= \langle y, x^*\rangle_r \langle x^*, y\rangle_r + \langle xb, xb\rangle_r\cr
&\leq \|x^*\|_r^2\langle y,y\rangle_r + \|x\|_r^2b^*b\cr
&\leq \|x\|_m^2(\langle y,y\rangle_r + b^*b).\cr}$$
This actually shows that $\|\psi(x)\| \leq \|x\|_m$, and the converse
inequality follows from the computations
$$\langle\psi(x)(\langle x,x\rangle_r \oplus 0),
\psi(x)(\langle x,x\rangle_r \oplus 0)\rangle_r
= \langle x,x\rangle_r^3$$
(hence $\|\psi(x)\| \geq \|x\|_r$) and
$$\langle \psi(x)(0 \oplus x^*), \psi(x)(0 \oplus x^*)\rangle_r
= \langle x^*, x^*\rangle_r^2$$
(hence $\|\psi(x)\| \geq \|x^*\|_r = \|x\|_l$).
Also, $\psi(x)$ is adjointable and in fact $\psi(x)^* = \psi(x^*)$
by the computations
$$\eqalign{\langle \psi(x)(b \oplus y), (c\oplus z)\rangle_r
&= \langle \langle x^*, y\rangle_r \oplus xb, c \oplus z\rangle_r\cr
&= \langle y, x^*\rangle_r c + \langle xb,z\rangle_r\cr}$$
and
$$\eqalign{\langle (b \oplus y), \psi(x^*)(c\oplus z)\rangle_r
&= \langle b \oplus y, \langle x,z\rangle_r \oplus x^*c\rangle_r\cr
&= b^*\langle x,z\rangle_r + \langle y, x^*c\rangle_r.\cr}$$
Finally, we note that $\phi(a)\psi(x) = \psi(ax)$ and $\psi(x)\phi(a) =
\psi(xa)$; these are trivially verified. We list the preceding
facts in the following theorem.
\bigskip

\noindent {\bf Theorem 12.} {\it Let $A$ be a C*-algebra and $E$ a Hilbert
$*$-bimodule over $A$. Then $\phi: A \to B(F)$ is an isometric
$*$-homomorphism, $\psi: E \to B(F)$ is an isometric linear embedding, and
for every $a \in A$ and $x, y \in E$ we have $\phi(a)\psi(x) = \psi(ax)$,
$\psi(x)\phi(a) = \psi(xa)$, and $\psi(x)^* = \psi(x^*)$.}\hfill\hal
\bigskip
\bigskip

Although in general $\phi(\la x,y\ra) \neq \psi(x)\psi(y)$, if $A$ has a
unit then we do have
$$\phi(\la x,y\ra)(1\oplus 0) = \la x,y\ra \oplus 0
= \psi(x)\psi(y)(1 \oplus 0)$$
for all $x,y \in E$.
\bigskip

Of course, there is an analogous left module version of this section's
construction.
\bigskip
\bigskip

\noindent {\bf 4. Duality.}
\bigskip

If $A = M$ is a von Neumann algebra, it is natural to focus attention on
Hilbert $*$-bimodules which are dual spaces. But given any Hilbert
$*$-bimodule $E$ over a von Neumann algebra, the linking
algebra construction can be modified by using the dual module $E_r'$ in
place of $E_r$; this has the consequence that $B(F)$ is a von Neumann
algebra [11], so that if $E$ is not a dual space we can
replace it with the weak* closure of $\psi(E)$ in $B(F)$. This is the
idea behind the main result of this section. We need some terminology first.
\bigskip

\noindent {\bf Definition 13.} Let $E$ be a Hilbert $*$-bimodule over a
von Neumann algebra $M$. We define the {\it $*$-weak topology} on $E$ to
be the weakest topology such that the maps $x \mapsto \la x,y\ra$ and
$x \mapsto \la y,x\ra$ are continuous from $E$ into $M$ for all $y \in E$.
\medskip

If the unit ball of $E$ is $*$-weakly compact, we say that $E$ is a
{\it dual} bimodule. If for any bounded, ultraweakly convergent net
$a_i \to a$ in $M$ and any $x \in E$ we have $a_ix \to ax$ $*$-weakly, then
we say that $E$ is {\it normal}. Finally, if $E$ is both normal and dual we
call it a {\it W* Hilbert $*$-bimodule}.\hfill\hal
\bigskip

It is standard that if $E$ is dual in the above sense then it is actually a
dual Banach space. Thus, given Proposition 5, the following proposition is
routinely verified.
\bigskip

\noindent {\bf Proposition 14.} {\it If $E$ is a dual Hilbert $*$-bimodule
over a von Neumann algebra then $B_l(E)$ and $B_r(E)$ are von Neumann
algebras. In either case
a bounded net of operators $(T_i)$ converges ultraweakly to
$T$ if and only if $T_ix \to Tx$ $*$-weakly for all $x \in E$.}\hfill\hal
\bigskip

The normality condition is symmetric, because $a_i x \to ax$ $*$-weakly
if and only if $x^*a_i^* \to x^*a^*$  $*$-weakly.
Also, note that part of normality is automatic: if $a_i \to a$ then
$$\la a_ix, y\ra = a_i\la x,y\ra \to a\la x, y\ra = \la ax, y\ra$$
for all $x,y \in E$. But $\la xa_i, y\ra \to \la xa, y\ra$ need not
always hold, for instance in the case of a module constructed as in
Example 9 using a non-normal conditional expectation of von Neumann algebras.
\medskip

If we modify the linking algebra construction by replacing $E_r$ with
$E_r'$ as suggested above, normality of $E$ is crucial because the
map $\phi: M \to B(F)$ is then ultraweakly continuous if and
only if $E$ is a normal module.
\medskip

We can now formalize the dualization procedure indicated at the start of
this section.
\bigskip

\noindent {\bf Lemma 15.} {\it Let $M$ and $N$ be von Neumann algebras and
let $A$ be an ultraweakly dense $*$-subalgebra of $M$. Suppose $\phi: A
\to N$ is an bounded $*$-homomorphism and $a_i \to 0$ boundedly and
ultraweakly in $A \subset M$ implies $\phi(a_i) \to 0$ ultraweakly in $N$.
Then $\phi$ extends to an ultraweakly continuous $*$-homomorphism from $M$
to $N$.}
\medskip

\noindent {\it Proof.} Let $M'$ be the ultraweak closure of
$A' = \{a \oplus \phi(a): a \in A\}$ in $M \oplus N$. Then the natural
projection $\pi_M: M' \to M$ hos no kernel and hence is an isomorphism,
and the map $\pi_N\circ \pi_M^{-1}: M \to N$ is ultraweakly continuous.
It is clear that the restriction of this map to $A$ agrees with
$\phi$.\hfill\hal
\bigskip

\noindent {\bf Theorem 16.} {\it Let $M$ be a von Neumann algebra, let $A$
be an ultraweakly dense $*$-subalgebra of $M$, and let $E$ be a
pre-Hilbert $*$-bimodule over $A$. Suppose that for any bounded net
$(a_i)$ in $A$ and any $x, y \in E$, $a_i \to 0$ ultraweakly implies
$\la xa_i, y\ra \to 0$ ultraweakly. Then $E$ modulo its null space densely
embeds in a unique W* Hilbert $*$-bimodule over $M$.
\medskip

In particular, if $E$ is a normal Hilbert $*$-bimodule over $M$ then it
densely embeds in a unique W* Hilbert $*$-bimodule over $M$.}
\medskip

\noindent {\it Proof.} First let $N = \{x \in E: \|x\|_m = 0\}$ and replace
$E$ with $E/N$. Now define $E_r'$ to be the set of right $A$-linear
maps from $E$ into $M$ which are bounded for the seminorm $\|\cdot\|_r$.
It follows from [11] that $E_r'$ is a self-dual right Hilbert module over
$M$ whose inner product extends $\la \cdot, \cdot\ra_r$ on $E_r$ when
$x \in E_r$ is identified with the map $y \mapsto \la x,y\ra_r$.
\medskip

Let $F = M \oplus E_r'$ be the direct sum of right Hilbert modules and
define maps $\phi: M \to B(F)$ and $\psi: E \to B(F)$ as in Section 3.
For a bounded net $(T_i)$ in $B(F)$, ultraweak convergence is equivalent to
ultraweak convergence of $\la T_i(x), y\ra_r$ in $M$ for all $x, y \in F$.
Thus Lemma 15 implies that there is an ultraweakly continuous extension
of $\phi|_A$ to $M$, and restriction of operators to $M \oplus 0$ shows
that this extension must be $\phi$. So $\phi$ is ultraweakly continuous.
\medskip

Define $E'$ to be the ultraweak closure of $\psi(E)$ in $B(F)$. This is
a bimodule over $M \cong \phi(M)$ via operator multiplication, and
normality and duality are trivial. It is also straightforward to check
that the bimodule structure of $E'$ extends that of $\psi(E) \cong E$.
The inner product and adjoint can either be extended from $E$ by
continuity or defined directly by $\la x,y\ra = x(y(1 \oplus 0))$
and operator adjoints.
\medskip

For uniqueness, let $E''$ be any other bimodule with the same properties,
and define a map $T: E' \to E''$ by $T(\lim_{E'} x_i) = \lim_{E''} x_i$
for any bounded universal net $(x_i)$ in $E$. This map is well-defined
and unitary since
$$\la \lim_{E''}x_i, y\ra = \lim_M \la x_i, y\ra = \la \lim_{E'}x_i, y\ra$$
for any $y \in E$, which is enough.\hfill\hal
\bigskip
\bigskip

\noindent {\bf 5. Locality.}
\bigskip

In this section we formulate a $*$-bimodule version of a condition on
Hilbert modules which was independently introduced in [17] and [19].
The purposes to which it was put in these two papers were very different,
and even the definitions are not obviously equivalent. (Their equivalence
follows from Proposition 5.4.2 of [17].) Our interest in the condition
is that it has strong consequences for the structure of the
bimodule which are analogous to facts about self-dual Hilbert modules [11].
\medskip

In [17] and [19] the property of interest was a C* version of the centered
condition given next. The appropriate $*$-bimodule property incorporates
a self-adjointness requirement, which we identified in the commutative case in
Theorem 11.
\bigskip

\noindent {\bf Definition 17.} Let $E$ be a W* Hilbert $*$-bimodule over a
von Neumann algebra $M$. The {\it center} of $E$ is the set
$$Z(E) = \{x \in E: ax = xa \hbox{ for all }a \in M\}.$$
We say that $E$ is {\it centered} if $MZ(E)$ is $*$-weakly dense in $E$.
We say that $E$ is {\it local} if it is centered and the inner product of
any two self-adjoint elements of $Z(E)$ is self-adjoint in $M$.\hfill\hal
\bigskip

In the next result we say that two subspaces $F, F' \subset E$ are
{\it orthogonal} if $\la x,y\ra = \la y,x\ra = 0$ for all $x \in F$ and
$y \in F'$. We also use the notation $Z(M)_{sa}$ or $Z(E)_{sa}$ for the set
of self-adjoint elements in $Z(M)$ or $Z(E)$.
\bigskip

\noindent {\bf Theorem 18.} {\it Let $E$ be a local W* Hilbert
$*$-bimodule over a von Neumann algebra $M$ and let $F$ be a centered,
$*$-weakly closed, self-adjoint sub-bimodule of $E$. Then there is another
centered, $*$-weakly closed, self-adjoint sub-bimodule $F'$ of $E$ which
is orthogonal to $F$ and such that $E = F + F'$.}
\medskip

\noindent {\it Proof.} Let $S$ be a subspace of $Z(F)_{sa}$ which is finitely
generated as a module over $Z(M)_{sa}$. We claim that we can find a finite
set $\{x_1, \ldots, x_n\}$ which spans $S$ over $Z(M)_{sa}$
such that $\la x_i, x_j\ra = 0$ whenever $i \neq j$. To see
this let $\{x_1, \ldots, x_n\}$ be any finite set which spans $S$ and
assume inductively that $\la x_i, x_j\ra = 0$ for $i,j \leq n - 1$,
$i \neq j$. Since
$$a\la x,y\ra = \la ax, y\ra = \la xa, y\ra = \la x, ay\ra =
\la x, ya\ra = \la x,y\ra a$$
for any $x,y \in Z(F)$ and $a \in M$, it follows that the inner product
of any two elements of $Z(F)$ is in $Z(M)$. Thus $Z(F)$ satisfies the
hypothesis of Theorem 11 as a Hilbert $*$-bimodule over $Z(M)$. Using the
conclusion of Theorem 11 it is easy to verify that
$$y = x_n - \sum_{i=1}^{n-1} {{\la x_n,x_i\ra}\over{\la x_i,x_i\ra}}x_i,$$
is well-defined and orthogonal to $x_i$ ($i \leq n - 1$) and
$\{x_1, \ldots, x_{n-1}, y\}$ spans $S$. This proves the claim.
\medskip

Now fix $y \in Z(E)_{sa}$. For any $S \subset Z(F)_{sa}$ as above, let
$\{x_1, \ldots, x_n\}$ verify the claim and define
$$y_S = \sum_{i = 1}^n {{\la y,x_i\ra}\over{\la x_i,x_i\ra}}x_i \in S.$$
This expression is sensible and $\|y_S\| \leq \|y\|$ by appeal to
Theorem 11. Direct the subspaces $S$ by inclusion and let $T(y)$ be a
cluster point of the net $(y_S)$; then $T(y) \in Z(F)_{sa}$ and
$\la y - T(y), x\ra = 0$ for all $x \in Z(F)_{sa}$, and as there is at most
one element of $Z(F)_{sa}$ which can have this property $T$ is a
well-defined (orthogonal) projection from $Z(E)_{sa}$ onto $Z(F)_{sa}$.
\medskip

Let $F'$ be the $*$-weak closure of the set $M\cdot{\rm ker}(T)$.
It is clear that $F'$ is orthogonal to
$F$. Also $Z(F)_{sa} + Z(F')_{sa} = Z(E)_{sa}$, so $F + F'$ is $*$-weakly
dense in $E$. So for any $x \in E$ we can find a bounded net
$(y_i + y_i')$ such that $y_i \in F$ and $y_i' \in F'$ and
$y_i + y_i' \to x$ $*$-weakly. By orthogonality, the nets $(y_i)$ and
$(y_i')$ are also bounded and so they have cluster points $y \in F$ and
$y' \in F'$, and $y + y' = x$ by continuity. So $E = F + F'$.\hfill\hal
\bigskip

As we mentioned earlier, there is a strong structure theorem for local
W* Hilbert $*$-bimodules. Let $\{p_i\}$ be a family of central
projections in a von Neumann algebra $M$. Then it is easy to check that the
algebraic direct sum of the family $\{p_iM\}$ is a pre-Hilbert
$*$-bimodule with respect to the bilinear form
$$\la \oplus a_i, \oplus b_i\ra = \sum a_ib_i$$
and involution
$$\big(\oplus a_i\big)^* = \oplus a_i^*$$
and that it satisfies the normality condition of Theorem 16. Thus it has
a W* Hilbert $*$-bimodule completion, which we denote $\bigoplus p_iM$; this
completion consists of those elements $\oplus a_i$ with the property that both
of the sums $\sum a_ia_i^*$ and $\sum a_i^*a_i$ converge ultraweakly. The
center of $\bigoplus p_iM$ is $\bigoplus Z(p_iM)$ and it is therefore
centered and local. The next theorem gives a converse to this fact.
\bigskip

\noindent {\bf Theorem 19.} {\it Let $E$ be a local W* Hilbert $*$-bimodule
over a von Neumann algebra $M$. Then there is a family $\{p_i\}$ of central
projections of $M$ such that $E \cong \bigoplus p_iM$.}
\medskip

\noindent {\it Proof.} As in the proof of Theorem 18, regard $Z(E)$ as a local
W* Hilbert $*$-module over $Z(M)$. Observe that for any $x \in Z(E)_{sa}$ the
sequence
$$(\la x, x\ra^{1/2} + n^{-1})^{-1} x$$
is bounded and if $y$ is a $*$-weak cluster point of this sequence then
$\la y,y\ra$ is a projection in $Z(M)$. Now let $\{x_i\}$ be a maximal family
of orthogonal elements of $Z(E)_{sa}$ with the property that $\la x_i, x_i\ra$
is a projection in $Z(M)$. It follows from Theorem 18 and the preceding
observation that
$M\cdot{\rm span}\{x_i\}$ is $*$-weakly dense in $E$. Also the sub-bimodules
$Mx_i$ are pairwise orthogonal. So it suffices to show that each $Mx_i$ is
isomorphic to $p_iM$ where $p_i = \la x_i, x_i\ra$. But the kernel of the
map $a \mapsto ax_i$ is an ultraweakly closed ideal of $M$, hence is of
the form $q_iM$ for some central projection $q_i$, and clearly $q_i = 1 - p_i$.
Also, for any $a, b \in p_iM$ we have $ab = \la ax_i, bx_i\ra$ and
$(ax_i)^* = a^*x_i$. So indeed $Mx_i \cong p_iM$ as Hilbert
$*$-bimodules.\hfill\hal
\bigskip
\bigskip

\noindent {\bf 6. Sauvageot's construction.}
\bigskip

Our most sophisticated class of examples of Hilbert $*$-bimodules, which were
the original motivation for this investigation, arise from a construction
in noncommutative geometry given in [17]. We now present a simplified
version of this construction which exhibits its symmetry and also shows
its resemblance to K\"ahler differentials (see e.g.\ [8]). Several
instances of the construction are detailed in [17].
\medskip

The ingredients of the construction are a von Neumann algebra $M$ and a
$C_0^*$-semigroup (see [3]) of completely positive maps $\phi_t: M \to M$
($t \geq 0$) such that
\smallskip

{\narrower{
\noindent (1) $\phi_0 = {\rm id}_M$,
\smallskip

\noindent (2) $\phi_t(1) = 1$ for all $t$, and
\smallskip

\noindent (3) the set
$$A_\infty = \{a \in M: a \in D(\Delta^n)\hbox{ for all }n\}$$
is an algebra, where $\Delta$ is the generator of $(\phi_t)$.
\smallskip}}

\noindent We also require the existence of a faithful, normal, semifinite
trace $\tau$ on $M$ such that
\smallskip

{\narrower{
\noindent (4) $\tau(a^*a) < \infty$ for all $a \in A_\infty$ and
\smallskip

\noindent (5) $\tau(a\phi_t(b)) = \tau(\phi_t(a)b)$ for all $a,b \in A_\infty$.
\medskip}}

Condition (5) is a noncommutative version of symmetry for Markov processes.
In [17] $\tau$ is only a weight, but in [18] it is also required to be a
trace.
\medskip

The construction proceeds as follows. For $a, b, c, d \in A_\infty$
define
$$(a\otimes b)^* = b^*\otimes a^*$$
and
$$\la a\otimes b, c\otimes d\ra = a\Delta(bc)d,$$
and extend both by linearity to $A_\infty\otimes A_\infty$. This bilinear
form does not satisfy $\la x, x^*\ra \geq 0$ on $A_\infty \otimes A_\infty$,
but it does hold on the sub-$A_\infty$-$A_\infty$-bimodule
$$E_0 = {\rm span}\{a \otimes bc - ab\otimes c: a, b, c \in A_\infty\}$$
([5], Theorem 14.7). Thus $E_0$ is a pre-Hilbert $*$-bimodule over $A_\infty$.
\medskip

The hypothesis of Theorem 16 is verified by observing that if $(a_i) \subset
A_\infty$ is bounded and $a_i \to 0$ ultraweakly in $M$ then
$$\eqalign{\tau(\la (a\otimes bc - &ab\otimes c)a_i,
(a'\otimes b'c' - a'b'\otimes c')\ra d)\cr
&= \tau(a\Delta(bca_ia')b'c'd - a\Delta(bca_ia'b')c'd
- ab\Delta(ca_ia')b'c'd + ab\Delta(ca_ia'b')c'd)\cr
&= \tau(bca_ia'\Delta(b'c'da) - bca_ia'b'\Delta(c'da)
- ca_ia'\Delta(b'c'dab) + ca_ia'b'\Delta(c'dab))\cr}$$
for all $a, b, c, a', b', c', d \in A_\infty$; the last expression converges
to zero since $\tau$ is normal, and therefore so does
$$\la (a\otimes bc - ab\otimes c)a_i, (a'\otimes b'c' - a'b'\otimes c')\ra$$
since $\tau$ is faithful and semifinite (see e.g.\ [9]). Thus $E_0$
densely embeds in a unique W* Hilbert $*$-bimodule $E$ over $M$.
\medskip

$E$ plays the role of the module of bounded measurable 1-forms, and we have
an exterior derivative $d_0: A_\infty \to E$ defined by
$$d_0(a) = i(1 \otimes a - a\otimes 1).$$
It is easy to check that $d_0$ is a self-adjoint derivation. In fact, it
is ultraweakly to $*$-weakly closable because $(a_i)$ and
$(d_0(a_i))$ bounded and $a_i \to 0$ ultraweakly imply that
$$\eqalign{\tau(\la i(1\otimes a_i - a_i\otimes 1),
&(a\otimes bc - ab\otimes c)\ra d)\cr
&= i\tau(\Delta(a_ia)bcd - \Delta(a_iab)cd
- a_i\Delta(a)bcd + a_i\Delta(ab)cd)\cr
&= i\tau(a_ia\Delta(bcd) - a_iab\Delta(cd)
- a_i\Delta(a)bcd + a_i\Delta(ab)cd)\cr}$$
converges to zero, and similarly for the inner product in reverse order,
which implies that $i(1\otimes a_i - a_i\otimes 1) \to 0$
$*$-weakly by the same reasoning as in the last paragraph. Thus, in the
terminology of [22] the closure $d$ of $d_0$ is a W*-derivation and its
domain is a noncommutative Lipschitz algebra.

%#
\bigskip
\bigskip

[1] C.\ A.\ Akemann, personal communication.
\medskip

[2] C.\ A.\ Akemann, J.\ Anderson, and G.\ K.\ Pedersen, Excising states
of C*-algebras, {\it Can.\ J.\ Math.\ \bf 38} (1986), 1239-1260.
\medskip

[3] O.\ Bratteli and D.\ W.\ Robinson, {\it Operator Algebras and
Quantum Statistical Mechanics I} (second edition), Springer-Verlag (1987).
\medskip

[4] A.\ Connes, C* alg\`ebres et g\'eom\'etrie diff\'erentielle, {\it
C.\ R.\ Acad.\ Sc.\ Paris (S\'er.\ A) \bf 290} (1980), A599-A604.
\medskip

[5] D.\ E.\ Evans and J.\ T.\ Lewis, Dilations of irreversible evolutions
in algebraic quantum theory, {\it Comm.\ Dublin Inst.\ Adv.\ Studies Ser.\ A
No.\ 24} (1977).
\medskip

[6] T.\ Falcone and M.\ Takesaki, Operator valued weights without
structure theory, {\it Trans.\ Amer.\ Math.\ Soc.\ \bf 351} (1999), 323-341.
\medskip

[7] J.\ M.\ G.\ Fell, {\it An extension of Mackey's method to Banach
$*$-algebraic bundles}, Mem.\ Amer.\ Math.\ Soc.\ 90 (1969).
\medskip

[8] R.\ Hartshorne, {\it Algebraic Geometry}, Springer-Verlag GTM {\bf 52}
(1977).
\medskip

[9] R.\ V.\ Kadison and J.\ R.\ Ringrose, {\it Fundamentals of the Theory
of Operator Algebras II}, Academic Press (1986).
\medskip

[10] C.\ Lance, {\it Hilbert C*-modules}, LMS Lecture Note Series {\bf 210},
Cambridge University Press (1995).
\medskip

[11] W.\ L.\ Paschke, Inner product modules over B*-algebras, {\it Trans.\
Amer.\ Math.\ Soc.\ \bf 182} (1973), 443-468.
\medskip

[12] G.\ K.\ Pedersen, {\it C*-algebras and their Automorphism Groups},
Academic Press (1979).
\medskip

[13] N.\ C.\ Phillips and N.\ Weaver, Modules with norms which take values
in a C*-algebra, {\it Pacific J.\ Math.\ \bf 185} (1998), 163-181
\medskip

[14] M.\ A.\ Rieffel, Induced representations of C*-algebras, {\it Adv.\
Math.\ \bf 13} (1974), 176-257.
\medskip

[15] ---------, Morita equivalence for C*-algebras and W*-algebras,
{\it J.\ Pure Appl.\ Algebra \bf 5} (1974), 51-96.
\medskip

[16] ---------, Morita equivalence for operator algebras, {\it Proc.\
Symp.\ Pure Math.\ \bf 38} (1982), 285-298.
\medskip

[17] J.-L.\ Sauvageot, Tangent bimodule and locality for dissipative
operators on C*-algebras, {\it Quantum Probability and App.\ IV}, Springer
LNM {\bf 1396} (1989), 322-338.
\medskip

[18] ---------, Quantum Dirichlet forms, differential calculus and
semigroups, {\it Quantum Probability and App.\ V}, Springer LNM
{\bf 1442} (1990), 334-346.
\medskip

[19] M.\ Skeide, Hilbert modules in quantum electro dynamics and
quantum probability, {\it Commun.\ Math.\ Phys.\ \bf 192} (1998), 569-604.
\medskip

[20] R.\ G.\ Swan, Vector bundles and projective modules, {\it Trans.\ Amer.\
Math.\ Soc.\ \bf 105} (1962), 264-277.
\medskip

[21] A.\ Takehashi, {\it Fields of Hilbert Modules}, Dissertation, Tulane
University (1971).
\medskip

[22] N.\ Weaver, Lipschitz algebras and derivations of von Neumann algebras,
{\it J.\ Funct.\ Anal.\ \bf 139} (1996), 261-300.
\medskip

[23] S.\ L.\ Woronowicz, Differential calculus on compact matrix pseudogroups
(quantum groups), {\it Comm.\ Math.\ Phys.\ \bf 122} (1989), 125-170.
\bigskip
\bigskip

\noindent Math Dept.

\noindent Washington University

\noindent St.\ Louis, MO 63130

\noindent nweaver@math.wustl.edu
\end